\documentclass{conm-p-l}

\numberwithin{equation}{section}

\usepackage{graphicx,amsthm}
\usepackage{verbatim,amsmath,amscd,amssymb,color}
\usepackage[mathscr]{eucal}
\input xy
\xyoption{all}
\begin{document}
\renewcommand{\labelenumi}{$($\roman{enumi}$)$}
\renewcommand{\labelenumii}{$(${\rm \alph{enumii}}$)$}
\font\germ=eufm10
\newcommand{\cI}{{\mathcal I}}
\newcommand{\cA}{{\mathcal A}}
\newcommand{\cB}{{\mathcal B}}
\newcommand{\cC}{{\mathcal C}}
\newcommand{\cD}{{\mathcal D}}
\newcommand{\cE}{{\mathcal E}}
\newcommand{\cF}{{\mathcal F}}
\newcommand{\cG}{{\mathcal G}}
\newcommand{\cH}{{\mathcal H}}
\newcommand{\cK}{{\mathcal K}}
\newcommand{\cL}{{\mathcal L}}
\newcommand{\cM}{{\mathcal M}}
\newcommand{\cN}{{\mathcal N}}
\newcommand{\cO}{{\mathcal O}}
\newcommand{\cR}{{\mathcal R}}
\newcommand{\cS}{{\mathcal S}}
\newcommand{\cT}{{\mathcal T}}
\newcommand{\cV}{{\mathcal V}}
\newcommand{\cX}{{\mathcal X}}
\newcommand{\cY}{{\mathcal Y}}
\newcommand{\fra}{\mathfrak a}
\newcommand{\frb}{\mathfrak b}
\newcommand{\frc}{\mathfrak c}
\newcommand{\frd}{\mathfrak d}
\newcommand{\fre}{\mathfrak e}
\newcommand{\frf}{\mathfrak f}
\newcommand{\frg}{\mathfrak g}
\newcommand{\frh}{\mathfrak h}
\newcommand{\fri}{\mathfrak i}
\newcommand{\frj}{\mathfrak j}
\newcommand{\frk}{\mathfrak k}
\newcommand{\frI}{\mathfrak I}
\newcommand{\fm}{\mathfrak m}
\newcommand{\frn}{\mathfrak n}
\newcommand{\frp}{\mathfrak p}
\newcommand{\fq}{\mathfrak q}
\newcommand{\frr}{\mathfrak r}
\newcommand{\frs}{\mathfrak s}
\newcommand{\frt}{\mathfrak t}
\newcommand{\fru}{\mathfrak u}
\newcommand{\frA}{\mathfrak A}
\newcommand{\frB}{\mathfrak B}
\newcommand{\frF}{\mathfrak F}
\newcommand{\frG}{\mathfrak G}
\newcommand{\frH}{\mathfrak H}
\newcommand{\frJ}{\mathfrak J}
\newcommand{\frN}{\mathfrak N}
\newcommand{\frP}{\mathfrak P}
\newcommand{\frT}{\mathfrak T}
\newcommand{\frU}{\mathfrak U}
\newcommand{\frV}{\mathfrak V}
\newcommand{\frX}{\mathfrak X}
\newcommand{\frY}{\mathfrak Y}
\newcommand{\fry}{\mathfrak y}
\newcommand{\frZ}{\mathfrak Z}
\newcommand{\frx}{\mathfrak x}
\newcommand{\rA}{\mathrm{A}}
\newcommand{\rC}{\mathrm{C}}
\newcommand{\rd}{\mathrm{d}}
\newcommand{\rB}{\mathrm{B}}
\newcommand{\rD}{\mathrm{D}}
\newcommand{\rE}{\mathrm{E}}
\newcommand{\rH}{\mathrm{H}}
\newcommand{\rK}{\mathrm{K}}
\newcommand{\rL}{\mathrm{L}}
\newcommand{\rM}{\mathrm{M}}
\newcommand{\rN}{\mathrm{N}}
\newcommand{\rR}{\mathrm{R}}
\newcommand{\rT}{\mathrm{T}}
\newcommand{\rZ}{\mathrm{Z}}
\newcommand{\bbA}{\mathbb A}
\newcommand{\bbB}{\mathbb B}
\newcommand{\bbC}{\mathbb C}
\newcommand{\bbG}{\mathbb G}
\newcommand{\bbF}{\mathbb F}
\newcommand{\bbH}{\mathbb H}
\newcommand{\bbP}{\mathbb P}
\newcommand{\bbN}{\mathbb N}
\newcommand{\bbQ}{\mathbb Q}
\newcommand{\bbR}{\mathbb R}
\newcommand{\bbV}{\mathbb V}
\newcommand{\bbZ}{\mathbb Z}
\newcommand{\adj}{\operatorname{adj}}
\newcommand{\Ad}{\mathrm{Ad}}
\newcommand{\Ann}{\mathrm{Ann}}
\newcommand{\rcris}{\mathrm{cris}}
\newcommand{\ch}{\mathrm{ch}}
\newcommand{\coker}{\mathrm{coker}}
\newcommand{\diag}{\mathrm{diag}}
\newcommand{\Diff}{\mathrm{Diff}}
\newcommand{\Dist}{\mathrm{Dist}}
\newcommand{\rDR}{\mathrm{DR}}
\newcommand{\ev}{\mathrm{ev}}
\newcommand{\Ext}{\mathrm{Ext}}
\newcommand{\cExt}{\mathcal{E}xt}
\newcommand{\fin}{\mathrm{fin}}
\newcommand{\Frac}{\mathrm{Frac}}
\newcommand{\GL}{\mathrm{GL}}
\newcommand{\Hom}{\mathrm{Hom}}
\newcommand{\hd}{\mathrm{hd}}
\newcommand{\rht}{\mathrm{ht}}
\newcommand{\id}{\mathrm{id}}
\newcommand{\im}{\mathrm{im}}
\newcommand{\inc}{\mathrm{inc}}
\newcommand{\ind}{\mathrm{ind}}
\newcommand{\coind}{\mathrm{coind}}
\newcommand{\Lie}{\mathrm{Lie}}
\newcommand{\Max}{\mathrm{Max}}
\newcommand{\mult}{\mathrm{mult}}
\newcommand{\op}{\mathrm{op}}
\newcommand{\ord}{\mathrm{ord}}
\newcommand{\pt}{\mathrm{pt}}
\newcommand{\qt}{\mathrm{qt}}
\newcommand{\rad}{\mathrm{rad}}
\newcommand{\res}{\mathrm{res}}
\newcommand{\rgt}{\mathrm{rgt}}
\newcommand{\rk}{\mathrm{rk}}
\newcommand{\SL}{\mathrm{SL}}
\newcommand{\soc}{\mathrm{soc}}
\newcommand{\Spec}{\mathrm{Spec}}
\newcommand{\St}{\mathrm{St}}
\newcommand{\supp}{\mathrm{supp}}
\newcommand{\Tor}{\mathrm{Tor}}
\newcommand{\Tr}{\mathrm{Tr}}
\newcommand{\wt}{\mathrm{wt}}
\newcommand{\Ab}{\mathbf{Ab}}
\newcommand{\Alg}{\mathbf{Alg}}
\newcommand{\Grp}{\mathbf{Grp}}
\newcommand{\Mod}{\mathbf{Mod}}
\newcommand{\Sch}{\mathbf{Sch}}\newcommand{\bfmod}{{\bf mod}}
\newcommand{\Qc}{\mathbf{Qc}}
\newcommand{\Rng}{\mathbf{Rng}}
\newcommand{\Top}{\mathbf{Top}}
\newcommand{\Var}{\mathbf{Var}}
\newcommand{\pmbx}{\pmb x}
\newcommand{\pmby}{\pmb y}
\newcommand{\gromega}{\langle\omega\rangle}
\newcommand{\lbr}{\begin{bmatrix}}
\newcommand{\rbr}{\end{bmatrix}}
\newcommand{\cd}{commutative diagram }
\newcommand{\SpS}{spectral sequence}
\newcommand\C{\mathbb C}
\newcommand\hh{{\hat{H}}}
\newcommand\eh{{\hat{E}}}
\newcommand\F{\mathbb F}
\newcommand\fh{{\hat{F}}}
\newcommand\Z{{\mathbb Z}}
\newcommand\Zn{\Z_{\geq0}}
\newcommand\et[1]{\tilde{e}_{#1}}
\newcommand\ft[1]{\tilde{f}_{#1}}

\def\ge{\frg}
\def\ji(#1,#2){j_{#1}^{(#2)}}
\def\AA{{\mathcal A}}
\def\al{\alpha}
\def\bq{B_q(\ge)}
\def\bqm{B_q^-(\ge)}
\def\bqz{B_q^0(\ge)}
\def\bqp{B_q^+(\ge)}
\def\beneme{\begin{enumerate}}
\def\beq{\begin{equation}}
\def\beqn{\begin{eqnarray}}
\def\beqnn{\begin{eqnarray*}}
\def\bfi{{\mathbf i}}
\def\bfii0{{\bf i_0}}
\def\bigsl{{\hbox{\fontD \char'54}}}
\def\bbra#1,#2,#3{\left\{\begin{array}{c}\hspace{-5pt}
#1;#2\\ \hspace{-5pt}#3\end{array}\hspace{-5pt}\right\}}
\def\cd{\cdots}
\def\ci(#1,#2){c_{#1}^{(#2)}}
\def\Ci(#1,#2){C_{#1}^{(#2)}}
\def\ch(#1,#2){c_{#2,#1}^{-h_{#1}}}
\def\cc(#1,#2){c_{#2,#1}}
\def\CC{\mathbb{C}}
\def\CBL{\cB_L(\TY(B,1,n+1))}
\def\CBM{\cB_M(\TY(B,1,n+1))}
\def\CVL{\cV_L(\TY(D,1,n+1))}
\def\CVM{\cV_M(\TY(D,1,n+1))}
\def\ddd{\hbox{\germ D}}
\def\del{\delta}
\def\Del{\Delta}
\def\Delr{\Delta^{(r)}}
\def\Dell{\Delta^{(l)}}
\def\Delb{\Delta^{(b)}}
\def\Deli{\Delta^{(i)}}
\def\Delre{\Delta^{\rm re}}
\def\ei{e_i}
\def\eit{\tilde{e}_i}
\def\eneme{\end{enumerate}}
\def\ep{\epsilon}
\def\eeq{\end{equation}}
\def\eeqn{\end{eqnarray}}
\def\eeqnn{\end{eqnarray*}}
\def\fit{\tilde{f}_i}
\def\FF{{\rm F}}
\def\ft{\tilde{f}_}
\def\gau#1,#2{\left[\begin{array}{c}\hspace{-5pt}#1\\
\hspace{-5pt}#2\end{array}\hspace{-5pt}\right]}
\def\gl{\hbox{\germ gl}}
\def\hom{{\hbox{Hom}}}
\def\ify{\infty}
\def\io{\iota}
\def\kp{k^{(+)}}
\def\km{k^{(-)}}
\def\llra{\relbar\joinrel\relbar\joinrel\relbar\joinrel\rightarrow}
\def\lan{\langle}
\def\lar{\longrightarrow}
\def\lm{\lambda}
\def\Lm{\Lambda}
\def\mapright#1{\smash{\mathop{\longrightarrow}\limits^{#1}}}
\def\Mapright#1{\smash{\mathop{\Longrightarrow}\limits^{#1}}}
\def\mm{{\bf{\rm m}}}
\def\mik(#1,#2){{}^{#1}{\bf m}_{#2}}
\def\mijk(#1,#2,#3){{}^{#1}{\bf m}^{(#2)}_{#3}}
\def\nd{\noindent}
\def\nn{\nonumber}
\def\nnn{\hbox{\germ n}}
\def\catob{{\mathcal O}(B)}
\def\cij(#1,#2){c_{#1}^{(#2)}}
\def\oint{{\mathcal O}_{\rm int}(\ge)}
\def\ot{\otimes}
\def\op{\oplus}
\def\opi{\ovl\pi_{\lm}}
\def\osigma{\ovl\sigma}
\def\ovl{\overline}
\def\plm{\Psi^{(\lm)}_{\io}}
\def\qq{\qquad}
\def\q{\quad}
\def\qed{\hfill\framebox[2mm]{}}
\def\QQ{\mathbb Q}
\def\qi{q_i}
\def\qii{q_i^{-1}}
\def\ra{\rightarrow}
\def\ran{\rangle}
\def\rlm{r_{\lm}}
\def\ssl{\hbox{\germ sl}}
\def\slh{\widehat{\ssl_2}}
\def\ti{t_i}
\def\tii{t_i^{-1}}
\def\til{\tilde}
\def\tm{\times}
\def\tri{\bigtriangleup}
\def\tt{\frt}
\def\TY(#1,#2,#3){#1^{(#2)}_{#3}}
\def\ua{U_{\AA}}
\def\ue{U_{\vep}}
\def\uq{U_q(\ge)}
\def\uqp{U'_q(\ge)}
\def\ufin{U^{\rm fin}_{\vep}}
\def\ufinp{(U^{\rm fin}_{\vep})^+}
\def\ufinm{(U^{\rm fin}_{\vep})^-}
\def\ufinz{(U^{\rm fin}_{\vep})^0}
\def\uqm{U^-_q(\ge)}
\def\uqmq{{U^-_q(\ge)}_{\bf Q}}
\def\uqpm{U^{\pm}_q(\ge)}
\def\uqq{U_{\bf Q}^-(\ge)}
\def\uqz{U^-_{\bf Z}(\ge)}
\def\ures{U^{\rm res}_{\AA}}
\def\urese{U^{\rm res}_{\vep}}
\def\uresez{U^{\rm res}_{\vep,\ZZ}}
\def\util{\widetilde\uq}
\def\uup{U^{\geq}}
\def\ulow{U^{\leq}}
\def\bup{B^{\geq}}
\def\blow{\ovl B^{\leq}}
\def\vep{\varepsilon}
\def\vp{\varphi}
\def\vpi{\varphi^{-1}}
\def\VV{{\mathcal V}}
\def\xii{\xi^{(i)}}
\def\xij(#1,#2){x_{#1}^{(#2)}}
\def\Xij(#1,#2){X_{#1}^{(#2)}}
\def\Xdij(#1,#2){{X'}_{#1}^{(#2)}}
\def\Xiioi{\Xi_{\io}^{(i)}}
\def\xxi(#1,#2,#3){\displaystyle {}^{#1}\Xi^{(#2)}_{#3}}
\def\xx(#1,#2){\displaystyle {}^{#1}\Xi_{#2}}
\def\W1{W(\varpi_1)}
\def\WW{{\mathcal W}}
\def\wt{{\rm wt}}
\def\wtil{\widetilde}
\def\what{\widehat}
\def\wpi{\widehat\pi_{\lm}}
\def\ZZ{\mathbb Z}
\def\RR{\mathbb R}

\def\m@th{\mathsurround=0pt}
\def\fsquare(#1,#2){
\hbox{\vrule$\hskip-0.4pt\vcenter to #1{\normalbaselines\m@th
\hrule\vfil\hbox to #1{\hfill$\scriptstyle #2$\hfill}\vfil\hrule}$\hskip-0.4pt
\vrule}}

\newtheorem{pro}{Proposition}[section]
\newtheorem{thm}[pro]{Theorem}
\newtheorem{lem}[pro]{Lemma}
\newtheorem{ex}[pro]{Example}
\newtheorem{cor}[pro]{Corollary}
\newtheorem{conj}[pro]{Conjecture}
\theoremstyle{definition}
\newtheorem{df}[pro]{Definition}

\newcommand{\cmt}{\marginpar}
\newcommand{\seteq}{\mathbin{:=}}
\newcommand{\cl}{\colon}
\newcommand{\be}{\begin{enumerate}}
\newcommand{\ee}{\end{enumerate}}
\newcommand{\bnum}{\be[{\rm (i)}]}
\newcommand{\enum}{\ee}
\newcommand{\ro}{{\rm(}}
\newcommand{\rf}{{\rm)}}
\newcommand{\set}[2]{\left\{#1\,\vert\,#2\right\}}
\newcommand{\sbigoplus}{{\mbox{\small{$\bigoplus$}}}}
\newcommand{\ba}{\begin{array}}
\newcommand{\ea}{\end{array}}
\newcommand{\on}{\operatorname}
\newcommand{\eq}{\begin{eqnarray}}
\newcommand{\eneq}{\end{eqnarray}}
\newcommand{\hs}{\hspace*}

\title[Decorated Geometric Crystals and Polyhedral Realizations of type $D_n$]
{Decorated Geometric Crystals and Polyhedral Realizations of type $D_n$}

\author{Toshiki N\textsc{akashima}}
\address{Department of Mathematics, 
Sophia University, Kioicho 7-1, Chiyoda-ku, Tokyo 102-8554,
Japan}
\email{toshiki@sophia.ac.jp}
\thanks{supported in part by JSPS Grants 
in Aid for Scientific Research $\sharp 22540031$.}

\date{}


\begin{abstract}
This is a continuation of \cite{N4,N5}.
We shall show that for type $D_n$ 
the realization of crystal bases 
obtained from the decorated geometric crystals in \cite{BK2} coincides
 with the polyhedral realizations of crystal bases. 
\end{abstract}

\maketitle
\renewcommand{\thesection}{\arabic{section}}
\section{Introduction}
\setcounter{equation}{0}
\renewcommand{\theequation}{\thesection.\arabic{equation}}

In \cite{N2,NZ} we introduced the polyhedral realizations of crystal
bases and therein we described their explicit forms for type $A_n$,
$A^{(1)}_n$ and arbitrary rank 2 Kac-Moody algebras.
In \cite{Ho}, Hoshino described them explicitly for all simple Lie
algebras. It is a kind of realizations of crystal bases that 
is presented as a convex polyhedral domain 
in an infinite/finite $\bbZ$-lattice defined by some
system of linear inequalities.

In \cite{BK2}, Berenstein and Kazhdan introduced decorated 
geometric crystals for reductive algebraic groups. 
Geometric crystals are  a kind of geometric analogue to the
Kashiwara's crystal bases (\cite{BK}). 
Let $I$ be a finite index set. Associated with 
a Cartan matrix $A=(a_{i,j})_{i,j\in I}$, 
define the decorated geometric crystal $\cX=(\chi,f)$, which 
is a pair of geometric crystal 
$\chi=(X,\{e_i\}_i,\{\gamma_i\}_i,\{\vep_i\}_i)$ 
and  a certain special rational function $f$ with the condition \eqref{f}.

If we apply the procedure called ``ultra-discretization''(UD) to ``positive 
geometric crystals'', 
then we obtain certain free-crystals for
the Langlands dual type (\cite{BK,N}).  As for a positive decorated
geometric crystal $(\chi,f,T',\theta)$
applying  UD to the function $f$ and
considering the convex polyhedral domain defined by the inequality 
$UD(f)\geq0$, we get the crystal with the property
``normal''(\cite{K3}) and furthermore as a connected component with 
the highest weight $\lm$, 
we obtain the Langlands dual Kashiwara's crystal $B(\lm)$.

This result makes us recall the polyhedral realization of crystal
bases  as introduced above since it has similar way to get the crystal 
$B(\lm)$ from certain free-crystals, defined by the system of linear
inequalities. 
Thus, one of the main aims of this article is to show that the crystals
obtained by UD from positive decorated geometric crystals and the
polyhedral realizations of crystals coincide with each other for type
$D_n$, which is a continuation of \cite{N4} for type $A_n$.

In \cite{N5}, the explicit feature of the decoration function $f_B$ on
some special geometric crystal $TB^-_{w_0}$ for 
all classical Lie algebras. In \cite{N4} for type $A_n$
the coincidence of polyhedral 
realizations of crystals and the realizations of crystals by the
decoration functions are presented. Here, in this article 
we shall do these for type $D_n$ by using the results in \cite{Ho}.
For type $B_n$ and $C_n$ we shall discuss the relations of two
realizations in the forthcoming paper.

The organization of the article is as follows:
in Sect.2, we review the theory of crystals and their polyhedral
realizations and give the explicit feature of polyhedral realization of
type $D_n$. 
In Sect.3, first we introduce the theory of decorated geometric 
crystals following \cite{BK2} and 
their positive structures and ultra-discretization.  
Next, we define the decoration function by using 
elementary characters and certain special positive decorated geometric
crystal on ${\bbB}_{w_0}=TB^-_{w_0}$. 
Finally, the ultra-discretization of $TB^-_{w_0}$ is
described explicitly. 
We describe the function $f_B$ exactly for type
$D_n$ in Sect.4 following \cite{N5}. In Sect.5, for  type $D_n$ 
we shall see the coincidence of the polyhedral
realization $\Sigma_{\io_0}[\lm]$ and the ultra-discretization 
$B_{f_B,\Theta^-_{\bfii0}}(\lm)$.

\renewcommand{\thesection}{\arabic{section}}
\section{Crystals and polyhedral realizations}
\setcounter{equation}{0}
\renewcommand{\theequation}{\thesection.\arabic{equation}}

\subsection{Notations}

Let $\ge$ be
a  semi-simple Lie algebra over $\bbQ$
with a Cartan subalgebra
$\tt$, a weight lattice $P \subset \tt^*$, the set of simple roots
$\{\al_i: i\in I\} \subset \tt^*$,
and the set of coroots $\{h_i: i\in I\} \subset \tt$,
where $I$ is a finite index set.
Let $\lan h,\lm\ran=\lm(h)$ be the pairing between $\tt$ and $\tt^*$,
and $(\al, \beta)$ be an inner product on
$\tt^*$ such that $(\al_i,\al_i)\in 2\bbZ_{\geq 0}$ and
$\lan h_i,\lm\ran={{2(\al_i,\lm)}\over{(\al_i,\al_i)}}$
for $\lm\in\tt^*$ and $A:=(\lan h_i,\al_j\ran)_{i,j}$ is the associated Cartan matrix.
Let $P^*=\{h\in \tt: \lan h,P\ran\subset\ZZ\}$ and
$P_+:=\{\lm\in P:\lan h_i,\lm\ran\in\ZZ_{\geq 0}\}$.
We call an element in $P_+$ a {\it dominant integral weight}.
The quantum algebra $\uq$
is an associative
$\QQ(q)$-algebra generated by the $e_i$, $f_i \,\, (i\in I)$,
and $q^h \,\, (h\in P^*)$
satisfying the usual relations.
The algebra $\uqm$ is the subalgebra of $\uq$ generated 
by the $f_i$ $(i\in I)$.

For the irreducible highest weight module of $\uq$
with the highest weight $\lm\in P_+$, we denote it by $V(\lm)$
and its {\it crystal base} we denote $(L(\lm),B(\lm))$.
Similarly, for the crystal base of the algebra $\uqm$ we denote 
$(L(\ify),B(\ify))$ (see \cite{K0,K1}).
Let $\pi_{\lm}:\uqm\longrightarrow V(\lm)\cong \uqm/{\sum_i\uqm
f_i^{1+\lan h_i,\lm\ran}}$ be the canonical projection and 
$\widehat \pi_{\lm}:L(\ify)/qL(\ify)\longrightarrow L(\lm)/qL(\lm)$
be the induced map from $\pi_{\lm}$. Here note that 
$\widehat \pi_{\lm}(B(\ify))=B(\lm)\sqcup\{0\}$.

Let {\it crystal } be a combinatorial object 
defined in \cite{K3}, see also \cite{N2,N4,NZ}.
In fact, $B(\ify)$ and $B(\lm)$ are the typical examples 
of crystals.

Let $B_1$ and $B_2$ be crystals.
A {\sl strict morphism } of crystals $\psi:B_1\lar B_2$
is a map $\psi:B_1\sqcup\{0\} \lar B_2\sqcup\{0\}$
satisfying the following conditions: $\psi(0)=0$, 
$wt(\psi(b)) = wt(b)$,
$\vep_i(\psi(b)) = \vep_i(b),$
$\vp_i(\psi(b)) = \vp_i(b)$
if $b\in B_1$ and $\psi(b)\in B_2,$
and the map $\psi: B_1\sqcup\{0\} \lar B_2\sqcup\{0\}$
commutes with all $\eit$ and $\fit$.
An injective strict morphism is called an {\it embedding }of crystals.

Crystals  have very nice properties for 
tensor operations. Indeed, if $(L_i,B_i)$ is a crystal base of 
$\uq$-module $M_i$ ($i=1,2$), $(L_1\ot_A L_2, B_1\ot B_2)$
is a crystal base of $M_1\ot_{\QQ(q)} M_2$ (\cite{K1}).

\subsection{Polyhedral Realization of Crystals}
\label{poly-uqm}
Let us recall the results in \cite{N2,NZ}.

Consider the infinite $\bbZ$-lattice
\begin{equation}
\ZZ^{\ify}
:=\{(\cd,x_k,\cd,x_2,x_1): x_k\in\ZZ
\,\,{\rm and}\,\,x_k=0\,\,{\rm for}\,\,k\gg 0\}.
\label{uni-cone}
\end{equation}
We fix an infinite sequence of indices
$\io=\cd,i_k,\cd,i_2,i_1$ from $I$ such that
\begin{equation}
{\hbox{
$i_k\ne i_{k+1}$ and $\sharp\{k: i_k=i\}=\ify$ for any $i\in I$.}}
\label{seq-con}
\end{equation}

Given $\io$, we can define a crystal structure
on $\ZZ^{\ify}$ and denote it by $\ZZ^{\ify}_{\io}$ 
(\cite[2.4]{NZ}).

\begin{pro}[\cite{K3}, See also \cite{NZ}]
\label{emb}
There is a unique strict embedding of crystals
$($called Kashiwara embedding$)$
\begin{equation}
\Psi_{\io}:B(\ify)\hookrightarrow \ZZ^{\ify}_{\geq 0}
\subset \ZZ^{\ify}_{\io},
\label{psi}
\end{equation}
such that
$\Psi_{\io} (u_{\ify}) = (\cd,0,\cd,0,0)$, where 
$u_{\ify}\in B(\ify)$ is the vector corresponding to $1\in \uqm$.
\end{pro}

In the rest of this section,
suppose $\lm\in P_+$. 
Let $R_\lm:=\{r_\lm\}$ be the crystal with the single element $r_\lm$
satisfying the condition $\wt(r_\lm)=\lm$, $\vep_i(r_\lm)=-\lan h_i,\lm\ran$
and $\vp_i(r_\lm)=0$.
Here we define the map
\begin{equation}
\Phi_{\lm}:(B(\ify)\ot R_{\lm})\sqcup\{0\}\longrightarrow B(\lm)\sqcup\{0\},
\label{philm}
\end{equation}
by $\Phi_{\lm}(0)=0$ and $\Phi_{\lm}(b\ot r_{\lm})=\wpi(b)$ for $b\in B(\ify)$.
We set
$$
\wtil B(\lm):=
\{b\ot r_{\lm}\in B(\ify)\ot R_{\lm}\,|\,\Phi_{\lm}(b\ot r_{\lm})\ne 0\}.
$$

\begin{thm}[\cite{N2}]
\label{ify-lm}
\begin{enumerate}
\item
The map $\Phi_{\lm}$ becomes a surjective strict morphism of crystals
$B(\ify)\ot R_{\lm}\longrightarrow B(\lm)$.
\item
$\wtil B(\lm)$ is a subcrystal of $B(\ify)\ot R_{\lm}$, 
and $\Phi_{\lm}$ induces the
isomorphism of crystals $\wtil B(\lm)\mapright{\sim} B(\lm)$.
\end{enumerate}
\end{thm}

By Theorem \ref{ify-lm}, we have the strict embedding of crystals
$\Omega_{\lm}:B(\lm)(\cong \wtil B(\lm))\hookrightarrow B(\ify)\ot R_{\lm}.$
Combining $\Omega_{\lm}$ and the
Kashiwara embedding $\Psi_{\io}$,
we obtain the following:
\begin{thm}[\cite{N2}]
\label{embedding}
There exists the unique  strict embedding of crystals
\begin{equation}
\Psi_{\io}^{(\lm)}:B(\lm)\stackrel{\Omega_{\lm}}{\hookrightarrow}
B(\ify)\ot R_{\lm}
\stackrel{\Psi_{\io}\ot {\rm id}}{\hookrightarrow}
\ZZ^{\ify}_{\io}[\lm],
\label{Psi-lm}
\end{equation}
such that $\Psi^{(\lm)}_{\io}(u_{\lm})=(\cd,0,0,0)\ot r_{\lm}$.
\end{thm}

Consider the infinite dimensional vector space
$$
\QQ^{\ify}:=\{{x}=
(\cd,x_k,\cd,x_2,x_1): x_k \in \QQ\,\,{\rm and }\,\,
x_k = 0\,\,{\rm for}\,\, k \gg 0\},
$$
and its dual space $(\QQ^{\ify})^*:={\rm Hom}(\QQ^{\ify},\QQ)$.
We will write a linear form $\vp \in (\QQ^{\ify})^*$ as
$\vp({x})=\sum_{k \geq 1} \vp_k x_k$ ($\vp_j\in \QQ$)
for $x\in \QQ^{\ify}$.

Let 
$S_k=S_{k,\io}$ on $(\QQ^{\ify})^*$ be  
the piecewise-linear operator as in \cite{NZ} and set
\begin{eqnarray}
\Xi_{\io} &:=  &\{S_{j_l}\cd S_{j_2}S_{j_1}x_{j_0}\,|\,
l\geq0,j_0,j_1,\cd,j_l\geq1\},
\label{Xi_io}\\
\Sigma_{\io} & := &
\{x\in \ZZ^{\ify}\subset \QQ^{\ify}\,|\,\vp(x)\geq0\,\,{\rm for}\,\,
{\rm any}\,\,\vp\in \Xi_{\io}\}.
\end{eqnarray}

For a fixed $\io=(i_k)$ and a positive integer $k$, let $k^{(-)}$ be the maximum
number $m$ such that $m<k$ and $i_k=i_m$ if it exists, and 0 otherwise.
We impose on $\io$ the following positivity assumption:
\begin{equation}
{\hbox{if $\km=0$ then $\vp_k\geq0$ for any 
$\vp(x)=\sum_k\vp_kx_k\in \Xi_{\io}$}}.
\label{posi}
\end{equation}
\begin{thm}[\cite{NZ}]
Let $\io$ be a sequence of indices satisfying $(\ref{seq-con})$ 
and (\ref{posi}). Then we have 
${\rm Im}(\Psi_{\io})(\cong B(\ify))=\Sigma_{\io}$.
\end{thm}

For every $k \geq 1$, let 
$\what S_k = \what S_{k,\io}$ be  the piece-wise linear
operator for a linear function 
$\vp(x)=c+\sum_{k\geq 1}\vp_kx_k$
$(c,\vp_k\in\QQ)$ on $\QQ^{\ify}$ as in \cite{N2}.
For the fixed sequence $\io=(i_k)$, 
in case $\km=0$ for $k\geq1$, there exists unique $i\in I$ such that $i_k=i$.
We denote such $k$ by $\io^{(i)}$, namely, $\io^{(i)}$ is the first
number $k$ such that $i_k=i$.
Here for $\lm\in P_+$ and $i\in I$ we set
$\lm^{(i)}(x):=
-\beta^{(-)}_{\io^{(i)}}(x)=\lan h_i,\lm\ran-\sum_{1\leq j<\io^{(i)}}
\lan h_i,\al_{i_j}\ran x_j-x_{\io^{(i)}}.
$

For $\io$ and a dominant integral weight $\lm$,
let $\Xi_{\io}[\lm]$ be the set of all linear functions
generated by $\what S_k=\what S_{k,\io}$ 
from the functions $x_j$ ($j\geq1$)
and $\lm^{(i)}$ ($i\in I$), namely,
\begin{equation}
\begin{array}{ll}
\Xi_{\io}[\lm]&:=\{\what S_{j_l}\cd\what S_{j_1}x_{j_0}\,
:\,l\geq0,\,j_0,\cd,j_l\geq1\}
\\
&\cup\{\what S_{j_k}\cd \what S_{j_1}\lm^{(i)}(x)\,
:\,k\geq0,\,i\in I,\,j_1,\cd,j_k\geq1\}.
\end{array}
\label{Xi}
\end{equation}
Now we set
\begin{equation}
\Sigma_{\io}[\lm]
:=\{x\in \ZZ^{\ify}_{\io}[\lm](\subset \QQ^{\ify})\,:\,
\vp(x)\geq 0\,\,{\rm for \,\,any }\,\,\vp\in \Xi_{\io}[\lm]\}.
\label{Sigma}
\end{equation}

For a sequence $\io$ and a dominant integral weight $\lm$, a pair
$(\io,\lm)$ is called {\it ample}
if $\Sigma_{\io}[\lm]\ni\vec 0=(\cd,0,0)$.

\subsection{\bf $D_n$-case}

Let us identify the index set $I$ with $[1,n] := \{1,2,\cd,n\}$ 
in the standard way; thus, the Cartan matrix 
$(a_{i,j}= \lan h_i,\al_j\ran )_{1 \leq i,j \leq n}$ of type $D_n$ 
is given by 
\begin{equation}
a_{i,i}=\begin{cases}2&\text{ if }i=j,\\
-1&\text{ if }|i-j|=1, (i,j)\ne(n-1,n),(n,n-1) \text{ or } 
(i,j)=(n-2,n),(n,n-2),\\
0&\text{ otherwise. }
\end{cases}
\end{equation}

Following the result in \cite{Ho}, 
we shall describe the polyhedral realization of type $D_n$.

As we have seen in \cite{N0},
in semi-simple setting, to describe $B(\lm)$ we only need the finite
rank $\bbZ$-lattice such as $\bbZ^{N}$ where $N$ is the length of the 
longest element in the corresponding Weyl group, that is, 
for a reduced longest word $\io=(i_N,\cd,i_2,i_1)$
the crystal $B(\lm)$ can be embedded in $\bbZ^N\ot R_\lm$.
Thus, for type $D_n$ we will take the sequence 
$\io_0=(n\cd 2\,\,1)^{n-1}$
($N=n(n-1)$).
Associating with this sequence, we take a doubly-indexed variables
$(x^{(j)}_i)_{(j,i)\in[1,n-1]\times [1,n]}$.
Here we define the following set of linear forms:
\begin{eqnarray}
&&\qq\q\Xi_k:=\{\xij(i,k)-\xij(i-1,k+1)|i=1,\cd,n-2\}
\cup\{\xij(n-1,k)+\xij(n,k)-\xij(n-2,k+1),
\xij(n-2,k+1)-\xij(n-1,k+1)-\xij(n,k+1)\}\label{xi1}\\
&&\cup\{\xij(n,k)-\xij(n-1,k+1),\xij(n-1,k)-\xij(n,k+1)\}
\cup\{\xij(n+k-j-1,j)-\xij(n+k-j,j)|j=k+2,\cd,n-1\}\q
(k=1,\cd,n-2),\nn \\
&&\q\qq\Xi_{n-1}:=\{\xij(n-1,n-1)\},\q\Xi_n:=\{\xij(n,n-1)\},\label{xin-1}\\
&&\q\qq \Xi'_k:=\{\xij(k-j,j)-\xij(k-j+1,j)|j=1,2,\cd,k\}\q(k=1,2,\cd,n-2),
\label{xi'}
\end{eqnarray}
where we understand $\xij(i,j)=0$ unless $(j,i)\in[1,n-1]\times [1,n]$.

To define the sets of linear forms $\Xi'_{n-1}$ and $\Xi'_n$, we prepare 
certain combinatorial objects ``admissible patterns''(\cite{Ho}):
\[
 M_n:=\{\mu=(\mu_1,\mu_2,\cd,\mu_l)|0<\mu_l<\mu_{l-1}<\cd<\mu_1<n.\},
\]
where we call an element $\mu$ in $M_n$ an admissible pattern.
\begin{ex} For $n=4$, we have 
$M_4=\{(3,2,1),(3,2),(3,1),(3),(2,1),(2),(1)\}.$
\end{ex}
The following is trivial.
\begin{lem}\label{2^n-1}
The number of elements in $M_n$ is $2^{n-1}-1$.
\end{lem}
Using these admissible patterns, we shall define the linear forms:
For an admissible pattern $\mu=(\mu_1,\cd,\mu_l)\in M_n$, define
\begin{eqnarray}
&&\vp_{\mu}(x):=\begin{cases}\displaystyle
\sum_{k=1}^l(\Xij(n-\mu_k-1,\mu_k+k-1)-\Xij(n-\mu_k,\mu_k+k-1))
&\text{ if }\mu_l=1,\\
\displaystyle
\sum_{k=1}^l(\Xij(n-\mu_k-1,\mu_k+k-1)-\Xij(n-\mu_k,\mu_k+k-1))+\Xij(n,l)
&\text{ if }\mu_l\geq 2,
\end{cases}\label{vpmu}\\
&&\vp'_{\mu}(x):=\begin{cases}\displaystyle
\sum_{k=1}^l(\Xdij(n-\mu_k-1,\mu_k+k-1)-\Xdij(n-\mu_k,\mu_k+k-1))
&\text{ if }\mu_l=1,\\
\displaystyle
\sum_{k=1}^l(\Xdij(n-\mu_k-1,\mu_k+k-1)-\Xdij(n-\mu_k,\mu_k+k-1))+\Xdij(n,l)
&\text{ if }\mu_l\geq 2,
\end{cases}\label{vppmu}\\
&&\Xi'_{n-1}:=\{\vp_{\mu}(x)|\mu\in M_n\},\qq
\Xi'_{n}:=\{\vp'_{\mu}(x)|\mu\in M_n\},
\end{eqnarray}
where we set
\[
 \Xij(i,j)=\begin{cases}
\xij(n,j)&\text{ if $j$ is even and }i=n-1,\\
\xij(n-1,j)&\text{ if $j$ is even and }i=n,\\
\xij(i,j)&\text{ otherwise,}
\end{cases}
\qq
\Xdij(i,j)=\begin{cases}
\xij(n,j)&\text{ if $j$ is odd and }i=n-1,\\
\xij(n-1,j)&\text{ if $j$ is odd and }i=n,\\
\xij(i,j)&\text{ otherwise.}
\end{cases}
\]
Now, set 
$\Xi_{\io_0}:=\bigcup_{k=1,n}\Xi_k$.
\begin{thm} [\cite{Ho}]
\label{D_n}
Let $\lm=\sum_{1\leq k\leq n}\lm_k\Lm_k$ $(\lm_k\in \ZZ_{\geq0})$ 
be a dominant integral weight. 
In the above notations, the image 
${\rm Im} \,(\Psi^{(\lm)}_{\io_0})=\Sigma_{\io_0}[\lm]\cong B(\lm)$ 
is the set of all integer families
 $(\xij(i,j))_{(j,i)\in[1,n-1]\times[1,n]}$ 
given by 
\begin{eqnarray}
&&\Sigma_{\io_0}[\lm]={\rm Im} \,(\Psi^{(\lm)}_{\io_0})\nn \\
&&\qq =\{x=(\xij(i,j))\,|\,f(x)\geq0 \text{ for any }f\in\Xi_{\io_0}\text{ and }
\lm_k+g(x)\geq 0\text{ for any }g\in\Xi'_k\,\,\,(k=1,\cd,n)\}.
\label{Blm-Dn}
\end{eqnarray}
\end{thm}

\renewcommand{\thesection}{\arabic{section}}
\section{Decorated geometric crystals}
\setcounter{equation}{0}
\renewcommand{\theequation}{\thesection.\arabic{equation}}

The basic reference for this section is \cite{BK,BK2,N,N4,N5}.
\subsection{Definitions}

Let  
 $A=(a_{ij})_{i,j\in I}$ be 
an indecomposable Cartan matrix
with a finite index set $I$.
Let $(\frt,\{\al_i\}_{i\in I},\{h_i\}_{i\in I})$ 
be the associated
root data 
satisfying $\al_j(h_i)=a_{ij}$.
Let $\ge=\ge(A)=\lan \frt,e_i,f_i(i\in I)\ran$ 
be the simple Lie algebra associated with $A$
over $\bbC$ and $\Del=\Del_+\sqcup\Del_-$
be the root system associated with $\ge$, where $\Del_{\pm}$ is 
the set of positive/negative roots.

Define the simple reflections $s_i\in{\rm Aut}(\frt)$ $(i\in I)$ by
$s_i(h)\seteq h-\al_i(h)h_i$, which generate the Weyl group $W$.
Let $G$ be the simply connected simple algebraic group 
over $\bbC$ whose Lie algebra is $\ge=\frn_+\oplus \frt\oplus \frn_-$.
Let $U_{\al}\seteq\exp\ge_{\al}$ $(\al\in \Del)$
be the one-parameter subgroup of $G$.
The group $G$ (resp. $U^\pm$) is generated by 
$\{U_{\al}|\al\in \Del\}$ 
(resp. $\{U_\al|\al\in\Del_{\pm}$).
Here $U^\pm$ is a unipotent radical of $G$ and 
${\rm Lie}(U^\pm)=\frn_{\pm}$.
For any $i\in I$, there exists 
a unique group homomorphism
$\phi_i\cl SL_2(\bbC)\rightarrow G$ such that
\[
\phi_i\left(
\left(
\begin{array}{cc}
1&t\\
0&1
\end{array}
\right)\right)=\exp(t e_i),\,\,
 \phi_i\left(
\left(
\begin{array}{cc}
1&0\\
t&1
\end{array}
\right)\right)=\exp(t f_i)\qquad(t\in\bbC).
\]
Set $\al^\vee_i(c)\seteq
\phi_i\left(\left(
\begin{smallmatrix}
c&0\\
0&c^{-1}\end{smallmatrix}\right)\right)$,
$x_i(t)\seteq\exp{(t e_i)}$, $y_i(t)\seteq\exp{(t f_i)}$, 
$G_i\seteq\phi_i(SL_2(\bbC))$,
$T_i\seteq \alpha_i^\vee(\bbC^\times)$ 
and 
$N_i\seteq N_{G_i}(T_i)$. Let
$T$ be a maximal torus of $G$ 
which has $P$ as its weight lattice and 
Lie$(T)=\frt$.
Let 
$B^{\pm}(\supset T)$ be the Borel subgroup of $G$.
We have the isomorphism
$\phi:W\mapright{\sim}N/T$ defined by $\phi(s_i)=N_iT/T$.
An element $\ovl s_i:=x_i(-1)y_i(1)x_i(-1)$ is in 
$N_G(T)$, which is a representative of 
$s_i\in W=N_G(T)/T$.

\begin{df}
\label{def-gc}
Let $X$ be an affine algebraic variety over $\bbC$, 
$\gamma_i$, $\vep_i, f$ 
$(i\in I)$ rational functions on $X$, and 
$e_i:\bbC^\times\times X\to X$ a unital rational $\bbC^\times$-action.
A 5-tuple $\chi=(X,\{e_i\}_{i\in I},\{\gamma_i,\}_{i\in I},
\{\vep_i\}_{i\in I},f)$ is a 
$G$ (or $\ge$)-{\it decorated geometric crystal} 
if
\begin{enumerate}
\item
$(\{1\}\times X)\cap dom(e_i)$ 
is open dense in $\{1\}\times X$ for any $i\in I$, where
$dom(e_i)$ is the domain of definition of
$e_i\cl\C^\times\times X\to X$.
\item
The rational functions  $\{\gamma_i\}_{i\in I}$ satisfy
$\gamma_j(e^c_i(x))=c^{a_{ij}}\gamma_j(x)$ for any $i,j\in I$.
\item
The function $f$ satisfies
\begin{equation}
f(e_i^c(x))=f(x)+{(c-1)\vp_i(x)}+{(c^{-1}-1)\vep_i(x)},
\label{f}
\end{equation}
for any $i\in I$ and $x\in X$,  where $\vp_i:=\vep_i\cdot\gamma_i$.
\item
$e_i$ and $e_j$ satisfy the following relations:
\[
 \begin{array}{lll}
&\hspace{-20pt}e^{c_1}_{i}e^{c_2}_{j}
=e^{c_2}_{j}e^{c_1}_{i}&
{\rm if }\,\,a_{ij}=a_{ji}=0,\\
&\hspace{-20pt} e^{c_1}_{i}e^{c_1c_2}_{j}e^{c_2}_{i}
=e^{c_2}_{j}e^{c_1c_2}_{i}e^{c_1}_{j}&
{\rm if }\,\,a_{ij}=a_{ji}=-1,\\
&\hspace{-20pt}
e^{c_1}_{i}e^{c^2_1c_2}_{j}e^{c_1c_2}_{i}e^{c_2}_{j}
=e^{c_2}_{j}e^{c_1c_2}_{i}e^{c^2_1c_2}_{j}e^{c_1}_{i}&
{\rm if }\,\,a_{ij}=-2,\,
a_{ji}=-1,\\
&\hspace{-20pt}
e^{c_1}_{i}e^{c^3_1c_2}_{j}e^{c^2_1c_2}_{i}
e^{c^3_1c^2_2}_{j}e^{c_1c_2}_{i}e^{c_2}_{j}
=e^{c_2}_{j}e^{c_1c_2}_{i}e^{c^3_1c^2_2}_{j}e^{c^2_1c_2}_{i}
e^{c^3_1c_2}_je^{c_1}_i&
{\rm if }\,\,a_{ij}=-3,\,
a_{ji}=-1.
\end{array}
\]
\item
The rational functions $\{\vep_i\}_{i\in I}$ satisfy
$\vep_i(e_i^c(x))=c^{-1}\vep_i(x)$ and 
$\vep_i(e_j^c(x))=\vep_i(x)$ if $a_{i,j}=a_{j,i}=0$.
\end{enumerate}
\end{df}
We call the function $f$ in (iii) the {\it decoration} of $\chi$ and 
the relations in (iv) are called 
{\it Verma relations}.
If $\chi=(X,\{e_i\},\,\{\gamma_i\},\{\vep_i\})$
satisfies the conditions (i), (ii), (iv) and (v), 
we call $\chi$ a {\it geometric crystal}.
{\sl Remark.}
The definitions of $\vep_i$ and $\vp_i$ are different from the ones in 
e.g., \cite{BK2} since we adopt the definitions following
\cite{KNO,KNO2}. Indeed, if we flip $\vep_i\to \vep_i^{-1}$ and 
$\vp_i\to \vp_i^{-1}$, they coincide with ours.

Let $\cT_+$ be the category whose objects are algebraic tori over
$\bbC$ and whose morphisms are positive rational maps.
Then, let $\mathcal UD$ be the functor:
\[
\begin{array}{cccc}
{\mathcal UD}:& \cT_+&\longrightarrow &{\mathfrak Set}\\
&T& \mapsto &X_*(T)\\
&(g:T\to T') &\mapsto &(\what g:X_*(T)\to X_*(T')),
\end{array}
\]
which is given as in \cite{KNO,KNO2,N,N4}.

For a split algebraic torus $T$ over $\bbC$, let us denote 
its lattice of (multiplicative )characters(resp. co-characters) by $X^*(T)$ 
(resp. $X_*(T)$). By the usual way, we identify $X^*(T)$
(resp. $X_*(T)$) with the weight lattice $P$ (resp. the dual weight
lattice $P^*$). 

Let $\theta:T'\rightarrow X$ be a positive structure on 
a decorated geometric crystal $\chi=(X,\{e_i\}_{i\in I},
\{{\rm wt}_i\}_{i\in I},
\{\vep_i\}_{i\in I},f)$.
Applying the functor ${\mathcal UD}$ 
to positive rational morphisms
$e_{i,\theta}:\bbC^\tm \tm T'\rightarrow T'$ and
$\gamma_i\circ \theta,\vep_i\circ\theta,f\circ\theta:T'\ra \bbC$
(the notations are
as above), we obtain
\begin{eqnarray*}
\til e_i&:=&{\mathcal UD}(e_{i,\theta}):
\ZZ\tm X_*(T') \rightarrow X_*(T')\\
{\rm wt}_i&:=&{\mathcal UD}(\gamma_i\circ\theta):
X_*(T')\rightarrow \bbZ,\\
\wtil\vep_i&:=&{\mathcal UD}(\vep_i\circ\theta):
X_*(T')\rightarrow \bbZ,\\
\wtil f&:=& {\mathcal UD}(f \circ\theta):X_*(T')\rightarrow \bbZ.
\end{eqnarray*}
Now, for given positive structure $\theta:T'\rightarrow X$
on a geometric crystal 
$\chi=(X,\{e_i\}_{i\in I},\{{\rm wt}_i\}_{i\in I},
\{\vep_i\}_{i\in I})$, we associate 
the quadruple $(X_*(T'),\{\til e_i\}_{i\in I},
\{{\rm wt}_i\}_{i\in I},\{\wtil\vep_i\}_{i\in I})$
with a free crystal structure (see \cite[2.2]{BK}) 
and denote it by ${\mathcal UD}_{\theta,T'}(\chi)$.
We have the following theorem:

\begin{thm}[\cite{BK,BK2,N}]
For any geometric crystal 
$\chi=(X,\{e_i\}_{i\in I},\{\gamma_i\}_{i\in I},
\{\vep_i\}_{i\in I})$ and positive structure
$\theta:T'\rightarrow X$, the associated crystal 
${\mathcal UD}_{\theta,T'}(\chi)=
(X_*(T'),\{\wtil e_i\}_{i\in I},\{{\rm wt}_i\}_{i\in I},
\{\wtil\vep_i\}_{i\in I})$ 
is a Langlands dual Kashiwara's crystal.
\end{thm}
{\sl Remark.}
The definition of $\wtil\vep_i$ is different from the one in 
\cite[6.1.]{BK2} since our definition of $\vep_i$ corresponds to 
$\vep_i^{-1}$ in \cite{BK2}.

Now, for a positive decorated geometric crystal 
$\cX=((X,\{e_i\}_{i\in I},\{\gamma_i\}_{i\in I},
\{\vep_i\}_{i\in I},f),\theta,T')$, set 
\begin{equation}
 \wtil B_{\wtil f}:=\{\wtil x\in X_*(T')|\wtil f(\wtil x)\geq0\},
\label{btil}
\end{equation}
where $X_*(T')$ is identified with $\bbZ^{\dim(T')}$. Define 
\begin{equation}
B_{f,\theta}:=(\wtil B_{\wtil f},\wtil e_i|_{\wtil B_{\wtil f}},
\wt_i|_{\wtil B_{\wtil f}},
\wtil\vep_i|_{\wtil B_{\wtil f}})_{i\in I}.
\label{btheta}
\end{equation}
\begin{pro}[\cite{BK2}]
For a positive decorated geometric crystal 
$\cX=((X,\{e_i\}_{i\in I},\{\gamma_i\}_{i\in I},
\{\vep_i\}_{i\in I},f),\theta,T')$, the
 quadruple $B_{f,\theta}$ in (\ref{btheta}) is a normal crystal.
\end{pro}

\subsection{Characters}

Let $\what U:={\rm Hom}(U,\bbC)$ be the set of additive characters 
of $U$.
The {\it elementary character }$\chi_i\in \what U$ and  
the {\it standard regular character} $\chi^{\rm st}\in \what U$ are  
defined by 
\[
\chi_i(x_j(c))=\del_{i,j}\cdot c \q(c\in \bbC,\,\, i\in I),\qq
\chi^{st}=\sum_{i\in I}\chi_i.
\]
Let us define an anti-automorphism $\eta:G\to G$ 
by 
\[
 \eta(x_i(c))=x_i(c),\q  \eta(y_i(c))=y_i(c),\q \eta(t)=t^{-1}\q
 (c\in\bbC,\,\, t\in T),
\]
which is called the {\it positive inverse}.

The rational function $f_B$ on $G$ is defined by 
\begin{equation}
f_B(g)=\chi^{st}(\pi^+(w_0^{-1}g))+\chi^{st}(\pi^+(w_0^{-1}\eta(g))),
\label{f_B}
\end{equation}
for $g\in B\ovl w_0 B$, where $\pi^+:B^-U\to U$ is the projection by 
$\pi^+(bu)=u$.


\subsection{Decorated geometric crystal on $\bbB_w$}
For a Weyl group element $w\in W$, define $B^-_w$ by 
$B^-_w:=B^-\cap U\ovl w U$
and set  $\bbB_w:=TB^-_w$. 
Let $\gamma_i:\bbB_w\to\bbC$ be the rational function defined by 
\begin{equation}
\gamma_i:\bbB_w\,\,\hookrightarrow \,\,\,
B^-\,\,\mapright{\sim}\,\,T\times U^-\,\,
\mapright{\rm proj}\,\,\, T\,\,\,\mapright{\al_i^\vee}\,\,\,\bbC.
\label{gammai}
\end{equation}

For any $i\in I$, there exists the natural projection 
$pr_i:B^-\to B^-\cap \phi(SL_2)$. Hence, 
for any $x\in \bbB_w$ there exists unique
       $v=\begin{pmatrix}b_{11}&0\\b_{21}&b_{22}\end{pmatrix}
\in SL_2$ such that 
$pr_i(x)=\phi_i(v)$. Using this fact, we define 
the rational function $\vep_i$ on $\bbB_w$ by 
\begin{equation}
\vep_i(x)=\frac{b_{22}}{b_{21}}\q(x\in\bbB_w).
\label{vepi}
\end{equation}
The rational $\bbC^\times$-action $e_i$ on $\bbB_w$ is defined by
\begin{equation}
e_i^c(x):=x_i\left((c-1)\vp_i(x)\right)\cdot x\cdot
x_i\left((c^{-1}-1)\vep_i(x)\right)\qq
(c\in\bbC^\times,\,\,x\in \bbB_w),
\label{ei-action}
\end{equation}
if $\vep_i(x)$ is well-defined, that is, $b_{21}\ne0$, 
and $e_i^c(x)=x$ if $b_{21}=0$.\\
{\sl Remark.} The definition (\ref{vepi}) is different from the one in 
\cite{BK2}. Indeed, if we take $\vep_i(x)=b_{21}/b_{22}$, 
then it coincides with
the one in \cite{BK2}.
\begin{pro}[\cite{BK2}]
For any $w\in W$,
the 5-tuple $\chi:=(\bbB_w,\{e_i\}_i,\{\gamma_i\}_i,\{\vep_i\}_i,f_B)$
is a decorated geometric crystal, where 
$f_B$ is in (\ref{f_B}), $\gamma_i$ is in (\ref{gammai}), $\vep_i$ is in 
(\ref{vepi}) and $e_i$ is in (\ref{ei-action}).
\end{pro}

\def\ld{\ldots}
Taking the  longest Weyl group element $w_0\in W$, let 
$\bfii0=i_1\ld i_N$ be one of its reduced expressions and 
define the positive structure on $B^-_{w_0}$ as 
$\Theta^-_\bfii0:(\bbC^\times)^N\longrightarrow B^-_{w_0}$ by 
\[
 \Theta^-_\bfii0(c_1,\cd,c_N):=\pmby_{i_1}(c_1)\cd \pmby_{i_N}(c_N),
\]
where $\pmby_i(c)=y_i(c)\al^\vee(c^{-1})$, 
which is different from $Y_i(c)$ in 
\cite{N,N2,KNO,KNO2}. Indeed, $Y_i(c)=\pmby_i(c^{-1})$.
We also define the positive structure on $\bbB_{w_0}$ as
$T\Theta^-_\bfii0:T\times(\bbC^\times)^N\,\,\longrightarrow\,\,\bbB_{w_0}$  
by $T\Theta^-_\bfii0(t,c_1,\cd,c_N)
=t\Theta^-_\bfii0(c_1,\cd,c_N)$.
The explicit geometric crystal
structure on $\bbB_{w_0}=TB^-_{w_0}$ is described in \cite{N4,N5}.

\subsection{Ultra-Discretization of $\bbB_{w_0}=TB^-_{w_0}$}
Applying the ultra-discretization functor to
$\bbB_{w_0}$,
we obtain the free crystal ${\mathcal UD}(\bbB_{w_0})=X_*(T)\times \bbZ^{N}$,
where $N$ is the length of the longest element $w_0$.
Then define the map $\til h:{\mathcal UD}(\bbB_{w_0})
=X_*(T)\times \bbZ^{N}\,\,\to\,\,
X_*(T)(=P^*)$ as the projection to the left component and set
$B_{w_0}(\lm^\vee):=\til h^{-1}(\lm^\vee)$ and 
$B_{f_B,\Theta^-_{\bfii0}}(\lm^\vee):=B_{w_0}(\lm^\vee)\cap 
B_{{f_B},\Theta^-_{\bfii0}}$
for $\lm^\vee\in X_*(T)=P^*$. Let  $P^*_+:=\{h\in P^*|\Lm_i(h)\geq0
\text{ for any }i\in I\}$ and for $\lm^\vee=\sum_i\lm_ih_i\in P^*_+$, we define
 $\lm=\sum_i\lm_i\Lm_i\in P_+$. Then, we have 
\begin{thm}[\cite{BK2}]
The set $B_{f_B,\Theta^-_{\bfii0}}(\lm^\vee)$ is non-empty if 
$\lm^\vee\in P^*_+$
 and in that case, $B_{f_B,\Theta^-_{\bfii0}}(\lm^\vee)$ is isomorphic to 
$B(\lm)^L$, which is the Langlands dual crystal associated with 
$\ge^L$. 
\end{thm}
\begin{thm}[\cite{N4}]
\label{ud-tb}
Let $\lm^\vee\in P^*_+$. 
The explicit crystal structure of $B_{f_B,\Theta^-_\bfii0}(\lm^\vee)$
is as follows:
For $x=(x_1,\cd,x_N)\in B_{f_B,\Theta^-_\bfii0}(\lm^\vee)\subset \bbZ^N$, we
 have
\begin{equation}
\eit^n(x)=\begin{cases}(x'_1,\cd,x'_N)&\text{if }
{\mathcal UD}(f_B)(x'_1,\cd,x'_N)\geq0,\\
\qq\q 0&\text{otherwise,}
\end{cases}
\end{equation}
where 
\begin{eqnarray}
&& x'_j=x_j+\min\left(\min_{1\leq m<j,i_m=i}(n+\sum_{k=1 }^m
 a_{i_k,i}x_k),
\min_{j\leq m\leq N,i_m=i}(\sum_{k=1 }^m a_{i_k,i}x_k)\right)
\label{eitn} \\
&&\qq -\min\left(\min_{1\leq m\leq j,i_m=i}(n+\sum_{k=1 }^m
 a_{i_k,i}x_k),
\min_{j<m\leq N,i_m=i}(\sum_{k=1 }^m a_{i_k,i}x_k)\right),
\nn
\\
&&\wt_i(x)=\lm_i-\sum_{k=1}^Na_{i_k,i}x_k,\qq
\vep_i(x)=\max_{1\leq m\leq
 N,i_m=i}(x_m+\sum_{k=m+1}^{N}a_{i_k,i}x_k),
\label{wt-ud}
\end{eqnarray}
where $x=(x_1,\cd,x_N)$ belongs to $B_{f_B,\Theta^-_\bfii0}(\lm^\vee)$ if and
 only if ${\mathcal UD}(f_B)(x)\geq0$.
\end{thm}
It follows immediately from (\ref{eitn}):
\begin{lem}
\label{lem-ef}
Set $X_m:=\sum_{k=1 }^m a_{i_k,i}x_k$, ${\mathcal X}^{(i)}:=\min\{X_m|1\leq
 m\leq N,i_m=i\}$ $(i\in I)$ and $M^{(i)}:=\{l|1\leq l\leq N,i_l=i,
 X_l={\mathcal X}^{(i)}\}$. 
Define $m_e:=\max(M^{(i)})$ and $m_f:=\min(M^{(i)})$:
for $x\in B_{f_B,\Theta^-_\bfii0}(\lm^\vee)$, we get 
\begin{eqnarray}
&&\eit(x)=\begin{cases}(x_1,\cd,x_{m_e}-1,\cd,x_N)
&\text{ if }{\mathcal UD}(f_B)(x_1,\cd,x_{m_e}-1,\cd,x_N)\geq0,\\
0&\text{otherwise,}
\end{cases}
\label{th-eaction}\\
&&\fit(x)=\begin{cases}(x_1,\cd,x_{m_f}+1,\cd,x_N)
&\text{ if }{\mathcal UD}(f_B)(x_1,\cd,x_{m_f}+1,\cd,x_N)\geq0,\\
0&\text{otherwise.}
\end{cases}
\label{th-faction}
\end{eqnarray}
\end{lem}
Finally, due to the results in Sect.2 and in this section, 
we obtain the following theorem 
\begin{thm}[\cite{N4}]
\label{coin-thm}
If we have $B_{f_B,\Theta^-_\bfii0}(\lm^\vee)
=\Sigma_{\bfii0^{-1}}[\lm]^L$ as a subset of $\bbZ^N$,
then they are isomorphic to each other as crystals, where ${}^L$ means the
Langlands dual crystal, that is, it is defined by the transposed Cartan 
matrix and $\bfii0^{-1}$ refers to $\bfii0$ in opposite order.
\end{thm}
\renewcommand{\thesection}{\arabic{section}}
\section{Explicit form of the decoration $f_B$ of type $D_n$}
\setcounter{equation}{0}
\renewcommand{\theequation}{\thesection.\arabic{equation}}

\subsection{Generalized Minors and the function $f_B$}

For this subsection, see \cite{BFZ,BZ,BZ2}.
Let $G$ be a simply connected simple algebraic group over $\bbC$ and 
$T\subset G$ a maximal torus. 
Let  $X^*(T):=\Hom(T,\bbC^\times)$ and $X_*(T):=\Hom(\bbC^\times,T)$ be
the lattice of characters and co-characters respectively.
We identify $P$ (resp. $P^*$) with $X^*(T)$ 
(resp. $X_*(T)$) as above.
\begin{df}
For $\mu\in P_+$, the
{\it principal minor} $\Del_\mu:G\to\bbC$ is defined as
\[
 \Del_\mu(u^-tu^+):=\mu(t)\q(u^\pm\in U^\pm,\,\,t\in T).
\]
Let $\gamma,\del\in P$ be extremal weights such that 
$\gamma=u\mu$ and $\del=v\mu$ for some $u,v\in W$. 
Then the {\it generalized minor} $\Del_{\gamma,\del}$ is defined
by
\[
 \Del_{\gamma,\del}(g):=\Del_\mu(\ovl u^{-1}g\ovl v)
\q(g\in G),
\]
which is a regular function on $G$.
\end{df}
\begin{pro}[\cite{BK2}]
The function $f_B$ in (\ref{f_B}) is described as:
\begin{equation}
 f_B(g)=\sum_i\frac{\Del_{w_0\Lm_i,s_i\Lm_i}(g)
+\Del_{w_0s_i\Lm_i,\Lm_i}(g)}
{\Del_{w_0\Lm_i,\Lm_i}(g)}.
\end{equation}
\end{pro}
Let ${\bf i}=i_1\cd i_N$ be a reduced word for the longest Weyl
group element $w_0$. 
For $t\Theta_{\bf i}^-(c)\in \bbB_{w_0}=T\cdot B^-_{w_0}$, 
we get the following formula.
\begin{equation}
 f_B(t\Theta_{\bfii0}^-(c))
=\sum_i\Del_{w_0\Lm_i,s_i\Lm_i}(\Theta_{\bfii0}^-(c))
+\al_i(t)\Del_{w_0s_i\Lm_i,\Lm_i}(\Theta_{\bfii0}^-(c)).
\label{fb-th}
\end{equation}

\subsection{Explicit form of  $f_B(t\Theta_{\bf i}^-(c))$ of type $D_n$}
\label{fb-Dn}

The results in this subsection are given in \cite{N5}.
Fix the cyclic reduced longest word
$\bfii0=(1 2 \cd n-1\,n)^{n-1}$.
\begin{thm}[\cite{N5}]\label{thm-Dn1}
For $k\in\{1,2,\cd,n\}$ and $c=(\ci(i,j))=(\ci(1,1),\ci(2,1),\cd,
\ci(n-1,n-1),\ci(n,n-1))\in(\bbC^\times)^{n(n-1)}$, we have
\begin{eqnarray*}
&&\Del_{w_0\Lm_k,s_k\Lm_k}(\Theta^-_\bfii0(c))\\
&&\hspace{-20pt}=
\ci(1,k)+\frac{\ci(2,k)}{\ci(1,k+1)}+\cd+
\frac{\ci(n-2,k)}{\ci(n-3,k+1)}
+\frac{\ci(n-1,k)\ci(n,k)}{\ci(n-2,k+1)}
+\frac{\ci(n-2,k+1)}{\ci(n-1,k+1)\ci(n,k+1)}+
\frac{\ci(n,k)}{\ci(n-1,k+1)}+\frac{\ci(n-1,k)}{\ci(n,k+1)}
+\frac{\ci(n-3,k+2)}{\ci(n-2,k+2)}+\cd+
\frac{\ci(k,n-1)}{\ci(k+1,n-1)}\\
&&\qq\qq (k=1,2,\cd,n-2),\\
&&\Del_{w_0\Lm_{n-1},s_{n-1}\Lm_{n-1}}(\Theta^-_\bfii0(c))=\ci(n-1,n-1),\qq
\Del_{w_0\Lm_{n},s_n\Lm_{n}}(\Theta^-_\bfii0(c))=\ci(n,n-1).
\end{eqnarray*}
\end{thm}

\begin{thm}[\cite{N5}]\label{thm-Dn2}
Let $k$ be in $\{1,2,\cd,n-2\}$.
Then we have
\begin{eqnarray}
\Del_{w_0s_k\Lm_k,\Lm_k}(\Theta^-_\bfii0(c))
=\frac{1}{\ci(1,k)}+\sum_{j=1}^{k-1}\frac{\ci(k-j,j)}{\ci(k-j+1,j)}.
\end{eqnarray}
\end{thm}
The cases $k=n-1,n$ will be presented in \ref{Dn-del-n}.

\subsection{$\Del_{w_0s_{n-1}\Lm_{n-1},\Lm_{n-1}}(\Theta^-_\bfii0(c))$ 
and $\Del_{w_0s_n\Lm_n,\Lm_n}(\Theta^-_\bfii0(c))$}\label{Dn-del-n}

To state the results for 
$\Del_{w_0s_{n-1}\Lm_n,\Lm_{n-1}}(\Theta^-_\bfii0(c))$ and 
$\Del_{w_0s_n\Lm_n,\Lm_n}(\Theta^-_\bfii0(c))$,
we need to prepare the set of triangles $\tri'_n$ for type $D_n$:
\begin{equation}
\tri'_n:=\{(\ji(k,l)|1\leq k\leq l<
n)|1\leq \ji(k,l+1)\leq \ji(k,l)<\ji(k+1,l+1)\leq n
\q(1\leq k\leq l<n-1)\}.
\end{equation}
We visualize a triangle $(\ji(k,l))$ in $\tri'_n$ as follows:
\[(\ji(k,l))=
 \begin{array}{c}
\ji(1,1)\\
\ji(2,2)\ji(1,2)\\
\ji(3,3)\ji(2,3)\ji(1,3)\\
\cd\cd\cd\cd\\
\ji(n-1,n-1)\cd\ji(2,n-1)\ji(1,n-1)
\end{array}
\]
\begin{lem}
For any $k\in\{1,2,\cd,n-1\}$ there exists a unique $j$
($1\leq j\leq k+1$) such that the $k$th row of 
a triangle $(\ji(k,l))$ in $\tri_n'$ is in the following form:
\begin{equation}
k\text{-th row}\q
(\ji(k,k),\ji(k-1,k),\cd,\ji(2,k),\ji(1,k))
=(k+1, k,k-1,\cd,j+1,j-1,j-2,\cd,2,1),
\end{equation}
that is, we have $\ji(m,k)=m$ for $m<j$ and 
$\ji(m,k)=m+1$ for $m\geq j$.
\end{lem}
For a triangle $\del=(\ji(k,l))\in\tri'_n$, we list $j$'s as in the lemma:
$s(\del):=(s_1,s_2,\cd,s_{n-1})$, 
which we call the {\it label} of a triangle $\del$.
\begin{lem}\label{Del-2^n-1}
$|\tri'_n|=2^{n-1}$.
\end{lem}
\begin{lem}\label{s4-d}
For any $\del\in\tri_n'$ 
let $s(\del):=(s_1,s_2,\cd,s_{n-1})$ be its label.
Then 
\begin{enumerate}
\item
The label $s(\del)$  satisfies 
$1\leq s_k\leq k+1$ and $s_{k+1}=s_k$ or $s_k+1$ for $k=1,\cd,n-1$.
\item
Each  $k$-th row of a triangle $\del$ is in 
one of the following
     {\rm I, II, III, IV:}
\begin{enumerate}
\item[{\rm I.}] $s_{k+1}=s_k+1$ and $s_k=s_{k-1}$.
\item[{\rm II.}] $s_{k+1}=s_k$ and $s_k=s_{k-1}$.
\item[{\rm III.}] $s_{k+1}=s_k+1$ and $s_k=s_{k-1}+1$.
\item[{\rm IV.}] $s_{k+1}=s_k$ and $s_k=s_{k-1}+1$.
\end{enumerate}
Here we suppose that $s_0=1$ and $s_{n}=s_{n-2}+1$, which means that 
the 1st row must be in {\rm I,II or IV} 
and the $n-1$-th row is in {\rm I or IV}.
\end{enumerate}
\end{lem}
Now, we associate a Laurant 
monomial $m(\del)$ in variables $(\ci(i,j))_{(j,i)\in[1,n-1]\times[1,n]}$ 
with a triangle $\del=(\ji(k,l))$ by the
following way.
\begin{enumerate}
\item Let $s=(s_1\cd,s_{n-1})$  be the label of $\del\in\tri'_n$.
\item Suppose  $i$-th row is in the form I. If $1\leq i\leq n-2$, 
 then associate $\ci(i,s_i)$. For $i=n-1$, 
\begin{enumerate}
\item if $n+ s_{n-1}$ is even, then associate $\ci(n-1,s_{n-1})$.
\item if $n+ s_{n-1}$ is odd, then associate $\ci(n,s_{n-1})$.
\end{enumerate}
\item Suppose  $i$-th row is in the form IV. If $1\leq i\leq n-2$, 
 then associate ${\ci(i,s_i)}^{-1}$. For $i=n-1$, 
\begin{enumerate}
\item if $n+ s_{n-1}$ is even, then associate ${\ci(n,s_{n-1})}^{-1}$.
\item if $n+ s_{n-1}$ is odd, then associate ${\ci(n-1,s_{n-1})}^{-1}$.
\end{enumerate}
\item If  $i$-th row is in the form II or III, then associate 1.
\item  Take the product of all monomials as above for 
$1\leq i <n$, then
we obtain the monomial $m(\del)$ associated with $\del$. 
\end{enumerate}
Here we define the involutions $\xi$ and $\ovl{\,\,}$ on Laurant monomials in 
$(\ci(i,j))_{(j,i)\in[1,n-1]\times[1,n]}$.
\begin{equation}
\begin{array}{l}
\xi:\ci(n-1,k)\mapsto \ci(n,k),\q \ci(n,k)\mapsto \ci(n-1,k),
\ci(j,k)\mapsto \ci(j,k) \q(j\ne n-1,n),\\
{\,}^-:\ci(i,j)\mapsto {\ci(i,n-j)}^{-1}.
\end{array}
\label{op-d}
\end{equation}
Let us denote the special triangle such that 
$\ji(k,l)=k+1$(resp. $\ji(k,l)=k$) for any $k,l$ by $\del_h$ 
(resp. $\del_l$).
Indeed, we have
\begin{equation}
m(\del_h)=\begin{cases}
\ci(n-1,1)&\text{ if $n$ is odd,}\\
\ci(n,1)&\text{ if $n$ is even,}
\end{cases},\qq
m(\del_l)={\ci(n,n)}^{-1}.
\label{m-del-d}
\end{equation}
Now, we present $\Del_{w_0s_{n-1}\Lm_{n-1},\Lm_{n-1}}(\Theta^-_\bfii0(c))$ 
and $\Del_{w_0s_n\Lm_n,\Lm_n}(\Theta^-_\bfii0(c))$ for type $D_n$:
\begin{thm}[\cite{N5}]\label{thm-Lm-n-d} 
For type $D_n$, we have the explicit forms:
\begin{eqnarray}
&&\Del_{w_0s_n\Lm_n,\Lm_n}(\Theta^-_\bfii0(c))=
\sum_{\del\in\tri'_n\setminus\{\del_l\}}\ovl{m(\del)},\\
&&\Del_{w_0s_{n-1}\Lm_{n-1},\Lm_{n-1}}(\Theta^-_\bfii0(c))=
\sum_{\del\in\tri'_n\setminus\{\del_l\}}\ovl{\xi\circ m(\del)}.
\end{eqnarray}
\end{thm}
\begin{ex} 
The set of triangles $\tri'_5$ is as follows:
\begin{eqnarray*}
&&\begin{array}{c}
1\\21\\321\\4321
\end{array}\q
\begin{array}{c}
1\\21\\321\\5321
\end{array}\q
\begin{array}{c}
1\\21\\421\\5321
\end{array}\q
\begin{array}{c}
1\\31\\421\\5321
\end{array}\q
\begin{array}{c}
2\\31\\421\\5321
\end{array}\q
\begin{array}{c}
2\\31\\421\\5421
\end{array}\q
\begin{array}{c}
1\\31\\421\\5421
\end{array}\q
\begin{array}{c}
1\\21\\421\\5421
\end{array}\\
&&\begin{array}{c}
1\\31\\431\\5421
\end{array}\q
\begin{array}{c}
2\\31\\431\\5421
\end{array}\q
\begin{array}{c}
2\\32\\431\\5421
\end{array}\q
\begin{array}{c}
1\\31\\431\\5431
\end{array}\q
\begin{array}{c}
2\\31\\431\\5431
\end{array}\q
\begin{array}{c}
2\\32\\431\\5431
\end{array}\q
\begin{array}{c}
2\\32\\432\\5431
\end{array}\q
\begin{array}{c}
2\\32\\432\\5432
\end{array}
\end{eqnarray*}
and their labels $s(\del)$ are
\begin{eqnarray*}
&&(2,3,4,5),\q (2,3,4,4),\q(2,3,3,4),\q(2,2,3,4),\q(1,2,3,4),\q
(1,2,3,3),\q(2,2,3,3),\q(2,3,3,3),\\
&&(2,2,2,3),\q (1,2,2,3),\q(1,1,2,3),\q(2,2,2,2),\q(1,2,2,2),\q
(1,1,2,2),\q(1,1,1,2),\q(1,1,1,1).
\end{eqnarray*}
Then, we have the corresponding monomials $\ovl{m(\del)}$:
\begin{eqnarray*}
&& \ci(5,0),\q
\frac{\ci(3,1)}{\ci(5,1)},\q \frac{\ci(2,2)\ci(4,1)}{\ci(3,2)},\q
\frac{\ci(1,3)\ci(4,1)}{\ci(2,3)},\q\frac{\ci(4,1)}{\ci(1,4)},\q
\frac{\ci(3,2)}{\ci(1,4)\ci(4,2)},\q\frac{\ci(1,3)\ci(3,2)}{\ci(2,3)\ci(4,2)},\q
\frac{\ci(2,2)}{\ci(4,2)},\\
&&\frac{\ci(1,3)\ci(5,2)}{\ci(3,3)},\q
\frac{\ci(2,3)\ci(5,2)}{\ci(1,4)\ci(3,3)},\q\frac{\ci(5,2)}{\ci(2,4)},\q
\frac{\ci(1,3)}{\ci(5,3)},\q \frac{\ci(2,3)}{\ci(1,4)\ci(5,3)},\q
\frac{\ci(3,3)}{\ci(2,4)\ci(5,3)},\q \frac{\ci(4,3)}{\ci(3,4)},\q
\frac{1}{\ci(4,4)}.
\end{eqnarray*}
Thus, we have 
\begin{eqnarray*}
&&\Del_{w_0s_n\Lm_n,\Lm_n}(\Theta^-_\bfii0(c))=
\frac{\ci(3,1)}{\ci(5,1)}+ \frac{\ci(2,2)\ci(4,1)}{\ci(3,2)}+
\frac{\ci(1,3)\ci(4,1)}{\ci(2,3)}+\frac{\ci(4,1)}{\ci(1,4)}+
\frac{\ci(3,2)}{\ci(1,4)\ci(4,2)}+\frac{\ci(1,3)\ci(3,2)}{\ci(2,3)\ci(4,2)}+
\frac{\ci(2,2)}{\ci(4,2)}\\
&&\qq\qq\qq\qq+\frac{\ci(1,3)\ci(5,2)}{\ci(3,3)}+
\frac{\ci(2,3)\ci(5,2)}{\ci(1,4)\ci(3,3)}+\frac{\ci(5,2)}{\ci(2,4)}+
\frac{\ci(1,3)}{\ci(5,3)}+ \frac{\ci(2,3)}{\ci(1,4)\ci(5,3)}+
\frac{\ci(3,3)}{\ci(2,4)\ci(5,3)}+ \frac{\ci(4,3)}{\ci(3,4)}
+\frac{1}{\ci(4,4)}.
\end{eqnarray*}
\end{ex}

\renewcommand{\thesection}{\arabic{section}}
\section{Ultra-Discretization and Polyhedral Realizations of type $D_n$}
\setcounter{equation}{0}
\renewcommand{\theequation}{\thesection.\arabic{equation}}

\subsection{Ultra-discretization of ${\mathcal UD}(f_B)(x)$}

Since in this section we treat the type $D_n$, we identify 
$P$ with $P^*$ by $\lm\leftrightarrow \lm^\vee$.
Let us describe the explicit form of $B_{f_B,\Theta^-_{\bfii0}}(\lm)$ for
type $D_n$ applying the result in Theorem \ref{ud-tb} and 
show the coincidence of the crystals $B_{f_B,\Theta^-_{\bfii0}}(\lm)$
and $\Sigma_{\io_0}[\lm]$ in Sect.2 using Theorem \ref{coin-thm}.

For type $D_n$ let $\bfii0$ be as in \ref{fb-Dn}. 
We shall see the explicit form of ${\mathcal UD}(f_B)(x)$. Indeed, 
due to (\ref{fb-th}), it is sufficient to know 
the forms of $\Del_{w_0\Lm_k,s_k\Lm_k}(\Theta_\bfii0^-(c))$ 
and $\Del_{w_0s_k\Lm_k,\Lm_k}(\Theta_\bfii0^-(c))$, which are given in 
Theorem \ref{thm-Dn1}, Theorem \ref{thm-Dn2} and 
Theorem \ref{thm-Lm-n-d}. Thus, we have
\[
 {\mathcal UD}(f_B)(t,x)=\min_{1\leq j\leq n}
({\mathcal UD}(\Del_{w_0\Lm_j,s_j\Lm_j}(\Theta_\bfii0^-))(x),
{\mathcal UD}(\al_j(t))
+{\mathcal UD}(\Del_{w_0s_j\Lm_j,\Lm_j}(\Theta_\bfii0^-))(x))
\]
and 
\begin{eqnarray*}
&&{\mathcal UD}(\Del_{w_0\Lm_k,s_k\Lm_k}(\Theta_\bfii0^-))(x)\\
&&=\min(\{\xij(i,k)-\xij(i-1,k+1)|i=1,\cd,n-2\}
\cup\{\xij(n-1,k)+\xij(n,k)-\xij(n-2,k+1),
\xij(n-2,k+1)-\xij(n-1,k+1)-\xij(n,k+1)\}\\
&&\hspace{-8pt}\cup\{\xij(n,k)-\xij(n-1,k+1),\xij(n-1,k)-\xij(n,k+1)\}
\cup\{\xij(n+k-j-1,j)-\xij(n+k-j,j)|k+2\leq j<n\})=
\min\Xi_k\q (1\leq k\leq n-2),\nn \\
&&{\mathcal UD}(\Del_{w_0\Lm_{n-1},s_{n-1}\Lm_{n-1}}(\Theta_\bfii0^-))(x)
=\xij(n-1,n-1)\in\Xi_{n-1},\qq
{\mathcal UD}(\Del_{w_0\Lm_{n},s_{n}\Lm_{n}}(\Theta_\bfii0^-))(x)
=\xij(n,n-1)\in\Xi_n.\\
&&{\mathcal UD}(\Del_{w_0s_k\Lm_k,\Lm_k}(\Theta_\bfii0^-))(x))
=\min\{\xij(k-j,j)-\xij(k-j+1,j)|j=1,2,\cd,k\}=\min\Xi'_k\q(1\leq k\leq n-2).
\end{eqnarray*}
where $\xij(i,j)={\mathcal UD}(\cij(i,j))$ and 
we understand $\xij(0,m)=0$.
Here ${\mathcal UD}(\Del_{w_0s_k\Lm_k,\Lm_k}(\Theta_\bfii0^-))(x))$
$(k=n-1,n)$ is given in the next subsection.

\subsection{Explicit correspondence for $k=n-1,n$}
Now, let us see 
${\mathcal UD}(\Del_{w_0s_k\Lm_k,\Lm_k}(\Theta_\bfii0^-))(x))=\min\Xi'_k$
for $k=n-1,n$. For the purpose, we define the map 
from the set of admissible patterns to the set of labels
of triangles.
Let $M_n$ be the set of admissible patterns as before and 
$S_n$ be the set of labels of $\tri'_n\setminus\{\del_l\}$:
\[
 S_n:=\{s=(s_1,\cd,s_{n-1})|1\leq s_1\leq \cd \leq s_{n-1}<n,
 |s_k-s_{k-1}|\leq 1\,(k=1,2,\cd,n-1)\},
\]
where we set $s_0=1$ and $s_n=s_{n-2}+1$.
Set $K_n:=(n,n-1,\cd,3,2)\in M_n$ and for $m\in[1,n-1]$
\[
 U_{m}:=\underbrace{(1,1,\cd,1,1)}_{m}.
\]
For $a=(a_1,a_2,\cd,a_l)\in\bbZ^l$, set
\[
 a^{-1}:=(a_l,\cd,a_2,a_1).
\]
Now, let us define the map $F:M_n\to S_n$ by 
\begin{equation}
F(\mu)=F(\mu_1,\cd,\mu_l):=\left(K_n-U_{\mu_1}-\cd-U_{\mu_l}\right)^{-1}.
\end{equation}
We shall see the well-definedness and bijectivity of this map.
\begin{ex}
For $n=5$, we have 
\begin{eqnarray*}
&&F(4,3,1)=((5,4,3,2)-(1,1,1,1)-(1,1,1)-(1))^{-1}=(2,2,1,1)^{-1}=(1,1,2,2),\\
&& F(3,2)=((5,4,3,2)-(1,1,1)-(1,1))^{-1}=(3,2,2,2)^{-1}=(2,2,2,3).
\end{eqnarray*}
\end{ex}

\begin{lem}
The map $F:M_n\to S_n$ is a bijective map.
\end{lem}
{\sl Proof.}
Note that $\sharp M_n=\sharp S_n=2^{n-1}-1$. The injectivity of $F$ is 
trivial if it is well-defined. Thus, its bijectivity is evident by these
facts.
So, let us show the well-definedness of the map $F$, that is, 
$F(\mu)$ is in $S_n$ for any $\mu\in M_n$.
Let us show this for $\mu=(m_1,\cd,m_l)$ by the induction on $l$.
For $l=1$, set $\mu=(m)\in M_n$ $(1\leq m\leq n-1)$. Then we have 
\[
F(\mu)=(2,3,\cd,n-m,\underbrace{n-m,n-m+1,\cd,n-1}_m), 
\]
which turns out
to be an element in $S_n$. 
For $\mu'=(\mu_2,\cd,\mu_l)$ with $1\leq \mu_l<\cd<\mu_2<n-1$, assume
that $F(\mu'):=(a_1,\cd,a_{n-1})$ is an element in $S_n$.
These $a_1,\cd,a_{n-1}$ satisfy that 
\begin{equation}
a_j=j+1\text{ for }j=1,2,\cd,n-\mu_2\text{ and }
 a_{n-\mu_2+1}=n-\mu_2+1=a_{n-\mu_2}.
\label{apro}
\end{equation}
Let $\mu_1$ be a positive integer with $\mu_2<\mu_1<n$ and 
$\mu:=(\mu_1,\mu_2,\cd,\mu_l)$, which is in $M_n$. 
By the definition of $F$, we obtain 
\[
F(\mu)=(a_1,\cd,a_{n-\mu_1},a_{n-\mu_1+1}-1,\cd,a_{n-1}-1).
\]
Since $(a_1,a_2,\cd,a_{n-1})$ is in $S_n$, then we easily find that 
$F(\mu)$ is also in $S_n$. 
As mentioned above, 
the bijectivity is evident from the injectivity of $F$ and 
the fact that $\sharp M_n=\sharp S_n$ by 
Lemma \ref{2^n-1} and Lemma \ref{Del-2^n-1}.
\qed

Let us see the explicit correspondence of the set of linear forms 
$\Xi_k$ and $\Xi'_k$ to the ultra-discretizations of
$\Del_{w_0\Lm_k,s_k\Lm_k}$
and $\Del_{w_0s_k\Lm_k,\Lm_k}$.
\begin{lem}\label{lem-mono-lin}
For $\mu\in M_n$, let $\vp'_\mu(x)$(resp.   $\vp_\mu(x)$) be the linear forms 
as in \eqref{vppmu} (resp. \eqref{vpmu}), and for $s\in S_n$
let us denote the corresponding triangle by $\del_s$ and 
let $m(\del_s)$ be the corresponding monomial as in Sect.4.
Then we have for any $\mu\in M_n$
\begin{eqnarray}
&&{\mathcal UD}(\ovl{m(\del_{F(\mu)})})(x)=\vp'_\mu(x),\label{mono-lin1}\\
&&{\mathcal UD}(\ovl{\xi\circ m(\del_{F(\mu)})})=\vp_\mu(x).\label{mono-lin2}
\end{eqnarray}
\end{lem}
{\sl Proof.}
Let us show \eqref{mono-lin1}. We shall see the conditions that 
each variable $\xij(i,j)$ appears in the left-hand side and 
the right-hand side of \eqref{mono-lin1}.
Note that each coefficient of variable $\xij(i,j)$ in both sides of
\eqref{mono-lin1} can be $\pm1$ or $0$.
First, let us see the right-hand side.
In the case $i\ne n-1,n$, due to the explicit description in 
\eqref{vppmu} for $\mu=(\mu_1,\cd,\mu_l)\in M_n$
we can see that $\xij(i,j)$ appears in $\vp'_\mu(x)$ iff 
there exists $k$ such that $j=\mu_k+k-1$, $i=n-\mu_k-1$ and 
$\mu_k+k-1\ne \mu_{k-1}+k-2$ or $n-\mu_k-1\ne n-\mu_{k-1}$, which is 
equivalent to the condition 
$\mu_k<\mu_{k-1}-1$.
Similarly, $-\xij(i,j)$ appears in $\vp'_\mu(x)$ iff 
there exists $k$ such that $j=\mu_k+k-1$, $i=n-\mu_k$ and 
$\mu_{k+1}<\mu_k-1$.
Considering the cases $i=n-1,n$ similarly, we obtain the following 
results for a admissible pattern $\mu=(\mu_1,\cd,\mu_l)$:
\begin{lem}
\label{lem-vp}
\begin{enumerate}
\item
$\xij(i,j)$ ($i\ne n-1,n$) appears in $\vp'_\mu(x)$ iff 
there exists $k$ such that $j=\mu_k+k-1$, $i=n-\mu_k-1$ and 
$\mu_k<\mu_{k-1}-1$.
\item
$-\xij(i,j)$ ($i\ne n-1,n$) appears in $\vp'_\mu(x)$ iff 
there exists $k$ such that $j=\mu_k+k-1$, $i=n-\mu_k$ and 
$\mu_{k+1}<\mu_k-1$.
\item
$\xij(n-1,j)$ appears in $\vp'_\mu(x)$ iff 
$\mu_l\geq 2$, $l=j$ and $j$ is odd.
\item
$\xij(n,j)$ appears in $\vp'_\mu(x)$ iff 
$\mu_l\geq 2$, $l=j$ and $j$ is even.
\item
$-\xij(n-1,j)$ appears in $\vp'_\mu(x)$ iff 
$\mu_l=1$, $\mu_l+l-1=j$ and $j$ is even.
\item
$-\xij(n,j)$ appears in $\vp'_\mu(x)$ iff 
$\mu_l=1$, $\mu_l+l-1=j$ and $j$ is odd.
\end{enumerate}
\end{lem}
Next, let us see the conditions that 
$\cij(i,j)$ appears in $\ovl{m(\del_s)}$.
By the recipe in Sect.4 to obtain $m(\del)$, we have the following
results for a label $s=(s_1,\cd,s_{n-1})\in S_n$:
\begin{lem}\label{lem-m}
\begin{enumerate}
\item
$\cij(i,j)$ ($i\ne n-1,n$) appears in $\ovl{m(\del_s)}$ iff
$j=n-s_i$, $s_i=s_{i-1}+1$ and $s_{i+1}=s_i$.
\item
${\cij(i,j)}^{-1}$ ($i\ne n-1,n$) appears in $\ovl{m(\del_s)}$ iff
$j=n-s_i$, $s_i=s_{i-1}$ and $s_{i+1}=s_i+1$.
\item
$\cij(n-1,j)$ appears in $\ovl{m(\del_s)}$ iff
$j=n-s_{n-1}$, $n+s_{n-1}$ is odd and 
$s_{n-1}=s_{n-2}+1$, $s_n=s_{n-1}$.
\item
$\cij(n,j)$ appears in $\ovl{m(\del_s)}$ iff
$j=n-s_{n-1}$, $n+s_{n-1}$ is even and 
$s_{n-1}=s_{n-2}+1$, $s_n=s_{n-1}$.
\item
${\cij(n-1,j)}^{-1}$ appears in $\ovl{m(\del_s)}$ iff
$j=n-s_{n-1}$, $n+s_{n-1}$ is even and 
$s_{n-1}=s_{n-2}$, $s_n=s_{n-1}+1$.
\item
${\cij(n,j)}^{-1}$ appears in $\ovl{m(\del_s)}$ iff
$j=n-s_{n-1}$, $n+s_{n-1}$ is odd and 
$s_{n-1}=s_{n-2}$, $s_n=s_{n-1}+1$.
\end{enumerate}
\end{lem}
Here note that we assumed that $s_n=s_{n-2}+1$ beforehand.
We find that each condition in both lemmas are equivalent 
through the map $F$.
Let us see that the conditions in Lemma \ref{lem-vp} (i) and Lemma
\ref{lem-m} (i) are equivalent via the map $F$.
Assume that $j=\mu_k+k-1$, $i=n-\mu_k-1$ and 
$\mu_k<\mu_{k-1}-1$. Since $i=n-\mu_k-1$, 
for $K_n^{-1}=(k_1,\cd,k_{n-1})$ we have $k_i=n-\mu_k$. 
And then for $s=(s_j)=F(\mu)$ we have 
$s_i=k_i-(k-1)=n-\mu_k-k+1$, which means $n-s_i=\mu_k-k+1=j$.
It follows from $\mu_k<\mu_{k-1}-1$ that $s_i=s_{i-1}+1$ and
$s_{i+1}=s_i$.
The inverse implication can be shown similarly.
The equivalence for the cases (ii)-(vi) are also shown similarly.
Thus, we obtain \eqref{mono-lin1} and by considering in the similar
manner we get \eqref{mono-lin2}. 
Then we completed the proof of Lemma \ref{lem-mono-lin}
\qed

Thus, by Lemma \ref{lem-mono-lin} we have
\begin{pro}\label{mono-lin-n-1n}
For $k=n-1,n$ we have the following:
\begin{eqnarray}
{\mathcal UD}(\Del_{w_0s_k\Lm_k,\Lm_k}(\Theta_\bfii0^-))(x))
=\min{\Xi'_{k}}.
\end{eqnarray}
\end{pro}

Thus, if we identify ${\mathcal UD}(\al_k(t))$ with $\lm_k$
($k=1,2,\cd,n$), then 
the condition ${\mathcal UD}(f_B)(\lm,x)\geq0$ 
in Theorem \ref{ud-tb}
is equivalent to the condition in \eqref{Blm-Dn} in Theorem \ref{D_n}.

Now, by Theorem \ref{coin-thm} we obtain the following theorem:
\begin{thm}
\label{thm-b}
In $D_n$ case, for any dominant integral weight $\lm$, 
there exists an isomorphism of crystals 
 $B_{f_B,\Theta^-_{\bfii0}}[\lm]\cong\Sigma_{\io_0}[\lm]$
where $\Sigma_{\io_0}[\lm]$ is as in Theorem \ref{D_n} and 
$\io_0=\bfii0^{-1}$.
\end{thm}

\bibliographystyle{amsalpha}

\end{document}